\documentclass[11pt]{article}
\usepackage{amsfonts}
\usepackage{latexsym}
\usepackage{amsmath}
\usepackage{amssymb}
\usepackage{amsthm}
\newtheorem{theorem}{Theorem}[section]
\newtheorem{lemma}[theorem]{Lemma}

\newtheorem{proposition}[theorem]{Proposition}
\newtheorem{corollary}[theorem]{Corollary}

\theoremstyle{definition}

\theoremstyle{remark}

\newtheorem*{note*}{Note}



\begin{document}

\title{\bf Relative entropies for  convex bodies
}

\author{Justin  Jenkinson  and Elisabeth M. Werner 
\thanks{Partially supported by an NSF grant, a FRG-NSF grant and  a BSF grant}}

\date{}

\maketitle

\begin{abstract}
We introduce a new class of (not necessarily convex) bodies 
and show, among other things,  that these bodies provide yet another link between convex geometric analysis and information theory.
Namely, they give  geometric interpretations of the relative entropy of the cone measures of a convex body and its polar and related  quantities.  
\par
Such interpretations were first given by Paouris and Werner  for  symmetric  convex bodies in  the context of the  $L_p$-centroid bodies.  There,  the relative entropies appear after performing   second order expansions of certain expressions.
Now, no symmetry assumptions  are needed. Moreover, using  the new bodies, already  first order expansions make 
the relative entropies  appear. Thus,  these bodies detect  ``faster"  details  of the  boundary of a convex body than the  $L_p$-centroid bodies.

\end{abstract}

      \vskip 2cm \noindent

\section{Introduction.}

It has been observed in recent years that there is a close connection between convex geometric analysis and information
theory. An example is the parallel between geometric inequalities for convex
bodies and inequalities for probability densities. For instance, the Brunn-Minkowski
inequality and the entropy power inequality follow both in a very similar way from the
sharp Young inequality (see. e.g., \cite{DCT}).
\par
Further connections between convexity and information theory were established by  Lutwak, Yang, and Zhang \cite{LYZ2000, LYZ2002/1, LYZ3}). They showed in \cite{LYZ2002/1}  that the
Cramer-Rao inequality corresponds to an inclusion of the Legendre ellipsoid and the
polar $L_2$-projection body. 
The latter is a basic notion from the  $L_p$-Brunn-Minkowski
theory. 
This theory evolved rapidly over the last years and
due to a number of highly influential works (see, e.g., 
\cite{Ga3}, \cite{GaZ}, \cite{GrZh},  \cite{Hab} - 
\cite{NPRZ}, \cite{RuZ}, 
\cite{Schu} - 
\cite{WY1},
\cite{Z3}), it is now  a central part of modern convex geometry. In fact, this  affine geometry of bodies pertains to some questions that had been considered Euclidean in nature. For example, the famous Busemann-Petty Problem (finally laid to rest in
\cite{Ga1, GaKoSch, RuZ, Z1, Z2}), was shown to be an affine problem with the introduction of intersection bodies by Lutwak in \cite{Lu1988}.

\par
Two fundamental notions within the
$L_p$-Brunn-Minkowski theory are $L_p$-affine surface areas,   introduced by Lutwak in \cite{Lu2}, 
 and $L_p$-centroid bodies introduced by Lutwak and Zhang in \cite{LZ}. See Section 3 for the definition of those quantities.
Based on these quantities, 
Paouris and Werner  \cite{PW1} established yet another  relation between affine convex
geometry and information theory. They proved that the exponential of the
relative entropy of the cone measures of a  {\em symmetric} convex body and its polar equals a limit of
normalized $L_p$-affine surface areas. 
Moreover, they introduce  a new affine invariant  quantity $\Omega_K$ (see also Section 3 for the definition).
\par
Here we introduce a new class of (not necessarily convex) bodies which we call {\em mean width bodies}. We describe some of their properties. For instance, we show that they are always star convex and that they provide geometric interpretations of  
 $L_p$-affine surface areas. 
Many such geometric interpretations have been given (see e.g. \cite{MW2, SW4, SW2004, W2, WY, WY1}). The twist  here is that these new geometric interpretations of affine invariants for {\em convex }bodies  are expressed in terms of {\em not necessarily convex bodies} (see  also \cite{WY1}).
\par
More importantly, these bodies provide yet another link between convex geometric analysis and information theory:
The main result of the paper  shows that these new bodies  give  geometric interpretations of both, the relative entropy of the cone measures of a {\em not necessarily symmetric } convex body and its polar and the quantity $\Omega_K$. 
Such interpretations were first given by Paouris and Werner \cite{PW1}  only for {\em symmetric} convex bodies in  the context of the  $L_p$-centroid bodies.  There the relative entropies appear after performing  a {\em second order expansion} of certain expressions.
The remarkable fact now  is that, using the mean width bodies, already a {\em first order expansion} makes them appear. Thus,  these new bodies detect  ``faster"  details  of the  boundary of a convex body than the  $L_p$-centroid bodies.

\vskip 2mm
\noindent
\subsection{Notation}
We work in ${\mathbb R}^n$, which is
equipped with a Euclidean structure $\langle\cdot ,\cdot\rangle $.
We denote by $\|\cdot \|_2$ the corresponding Euclidean norm. $B^n_2(x, r)$ is the Euclidean ball centered at $x$ with radius $r$. 
We write $B_2^n=B^n_2(0, 1)$ for the Euclidean unit ball centered at $0$ and $S^{n-1}$ for the
unit sphere. Volume is denoted by $|\cdot |$. 
Throughout the paper, we will  assume that  the centroid of a
convex body $K$ in $\mathbb R^n$ is at the origin. 
$K^\circ=\{ y\in \mathbb{R}^n: \langle x, y \rangle \leq 1 \  \text{for all }  \  x \in K\}$ is the polar body of $K$.
\par
We  write $K\in
C^2_+$,  if $K $ has $C^2$ boundary $\partial K$  with everywhere strictly
positive Gaussian curvature $\kappa_K$. For a point $x \in \partial K$, the boundary of $K$, $N_K(x)$ is the outer unit normal 
in $x$ to $K$. $\mu_K$ is the usual surface area measure on $\partial K$.
 $\omega $ is the usual surface area measure on $S^{n-1}$
and  $\sigma $ its normalization: $\sigma (A) = \frac{\omega (A) }{\omega (S^{n-1})}$ for all Borel measurable sets $A \subset S^{n-1}$. 
\par
For  $u $ and $x$ in $ \mathbb{R}^n$, 
$H=H(x, \xi)$ is the hyperplane through $x$ orthogonal to $\xi$. $H^+=H^+(x, \xi)=\{y\in \mathbb{R}^n: \langle y, \xi \rangle \geq \langle x, \xi \rangle \}$ and  $H^-=H^-(x, \xi)=\{y\in \mathbb{R}^n: \langle y, \xi \rangle \leq \langle x, \xi \rangle \}$ are the two closed half spaces generated by 
$H$.
\par
Let  $K$ be a convex body in $\mathbb{R}^n$ and let  $u \in S^{n-1}$. Then
$h_{K}(u)$ is the support function of direction $u\in S^{n-1}$,
and $f_{K}(u)$ is the curvature function, i.e. the reciprocal of
the Gaussian curvature $\kappa _K(x)$ at this point $x \in
\partial K$ that has $u$ as outer normal.

\vskip 5mm
\section{Mean width bodies.}
The width $W(K)$ of a convex body $K$ in $\mathbb{R}^n$ is defined as
$$
W(K)= 2 \int_{S^{n-1}} h_K(u) d\sigma(u).
$$
Let $M$ and $K$ be convex bodies such that $0$ is the center of gravity of $K$ and $K \subset M$.
It is easy to see  \cite{GG}) that
\begin{equation}\label{w-MK}
 W(M)-W(K)=\frac{2}{\omega(S^{n-1})}\ \int_{K^\circ
\setminus M^\circ } \|\xi\|^{-(n+1)} d\xi.
\end{equation}
\vskip 3mm 
\noindent 
Let $f: K^\circ \rightarrow \mathbb{R}$ be a positive,  integrable function.
We generalize
(\ref{w-MK}) to
\begin{equation}\label{fw-MK}
 W_f(M)-W_f(K)= \frac{2}{\omega(S^{n-1})}\ \int_{K^\circ \setminus
M^\circ } f(\xi) d \xi
\end{equation}

\vskip 4mm 
\noindent
For the following easy lemma we will need another notation.
\par
\noindent
Let $\alpha \in \mathbb{R}$, $\alpha  \neq 0$. Let $f:S^{n-1} \rightarrow \mathbb {R}$ be a positive function.
Recall that $f$ is said to be  homogeneous of degree $  \alpha$,  if for all $r \geq 0$, 
$$
f(r u) =r^\alpha  f(u).
$$
\vskip 2mm 
\noindent
\begin{lemma} \label{easylemma}
Let $K$ and $M$ be  convex bodies in $\mathbb{R}^n$ such that
$0$ is the center of gravity of $K$ and $K \subset M$. Let $f: S^{n-1} \rightarrow \mathbb{R}$ be a positive, integrable 
function that is homogeneous of degree $ \alpha$.
\vskip 2mm
\noindent
(i) Let $\alpha \neq - n$. Then
$$
 W_f(M)-W_f(K)= \frac{2}{(\alpha+n)}\ \int_{S^{n-1}} f(u)
\left[\frac{1}{h_{K}^{\alpha+n}(u)} - \frac{1}{h_M^{\alpha+n}(u)}\right] d\sigma(u).
$$
\vskip 2mm
\noindent
(ii) Let $\alpha = - n$. Then
$$
 W_f(M)-W_f(K)= 2 \int_{S^{n-1}} f(u)
\log \left[\frac{h_{M}(u)}{ h_K(u)}\right] d\sigma(u).
$$
\end{lemma}
\vskip 3mm
\noindent
{\bf Proof.}
We use $\alpha$-homogeneity and get 
\begin{eqnarray*}
 W_f(M)-W_f(K) &=&  \frac{2}{\omega(S^{n-1})}\ \int_{K^\circ \setminus
M^\circ } f(\xi) d \xi \\
&=& \frac{2}{\omega(S^{n-1})}\ \int_{S^{n-1}} \int_{\frac{1}{h_{M}(u)}}
^{\frac{1}{h_{K}(u)}} f(ru) r^{n-1} dr d\omega (u)\\
& = &\frac{2}{\omega(S^{n-1})}\ \int_{S^{n-1}} \int_{\frac{1}{h_{M}(u)}} ^{\frac{1}{h_{K}(u)}} f(u) r^{n + \alpha -1} dr d\omega (u)
\end{eqnarray*}
Integration then yields (i) and (ii).
\vskip 4mm

Let $(X, \mu)$ be a measure space  and let  $dP=pd\mu$ and  $dQ=qd\mu$ be probability measures on $X$ that are  absolutely continuous with respect to the measure $\mu$. 
The {\em Kullback-Leibler divergence} or {\em relative entropy} from $P$ to $Q$ is defined as (see \cite{CT})
\begin{equation}\label{relent}
 D_{KL}(P\|Q)= \int_{X} p\log{\frac{p}{q}} d\mu.
\end{equation}
\vskip 4mm
If we let  $f(u) = \frac{1}{h_K^n(u)}$ (or $f(u) = \frac{1}{h_M^n(u)}$) 
in  Lemma \ref{easylemma} (ii),  then $f(r u)= \frac{r^{-n}}{h_K^n(u)}= r^{-n} f(u)$.
Thus  this $f$ is homogeneous of degree $-n$.  
\par
Let now $(X, \mu)=(S^{n-1}, \omega)$ and for convex bodies $K$ and $M$ in $\mathbb{R}^n$ put
\begin{equation}\label{pq}
p_K= \frac{1}{n |K^\circ| h_K^n}, \hskip 5mm p_M= \frac{1}{n |M^\circ | h_{M}^n}.
\end{equation}
Then $ dP_K=p_K  d\omega$ and $dP_M =  p_M d\omega $ are probability measures on $S^{n-1}$
and  Lemma \ref{easylemma} (ii) becomes
\begin{eqnarray*} \label{ent}
 W_{\frac{1}{h_K^n}}(M)-W_{\frac{1}{h_K^n}}(K)  &=&  \frac{2}{n } |K^\circ|\int_{S^{n-1}}  \frac{1}{|K^\circ | h_K^n} \log  \left(\frac{h_{M}^n}{ h_K^n}\right) d\sigma \nonumber \\
&=&\frac{2|K^\circ|}{ \omega(S^{n-1})} \int_{S^{n-1}} p_K\left(  \log \frac{p_K}{p_M} +\log \left( \frac{|K^\circ|}{|M^\circ|} \right)\right)d \omega \nonumber \\
&=&\frac{2|K^\circ|}{ \omega(S^{n-1})}\left( D_{KL}(P_K\|P_M) +\log \left( \frac{|K^\circ|}{|M^\circ|} \right) \right).
\end{eqnarray*}
Hence we get

\begin{corollary}
Let $K$ and $M$ be  convex bodies in $\mathbb{R}^n$ such that $K \subset M$ and let $p_K$ and $p_M$ be the probability densities given in (\ref{pq}). Then
$$
 \int_{K^\circ \setminus
M^\circ } \frac{1}{h_K^n(\xi)}  \frac{d \xi}{|K^\circ|} =D_{KL}(P_K\|P_M) + \log \left( \frac{|K^\circ|}{|M^\circ|} \right)
$$
\end{corollary}
\vskip 4mm 
We now want to apply the above considerations for a specific $M$. Namely, 
for $x \in \mathbb{R}^n$, let $K_x=[x,K]$ be the convex hull of
$x$ and $K$. For $x \in K$, $K_x=K$. Therefore, we will consider only  $x \notin K$.   Let $t \geq 0$ and let 
\begin{equation}\label{kt}
K[t]=\{ x\in \mathbb{R}^n: w(x) \leq t\}
\end{equation}
where
\begin{equation}\label{wx}
w(x) = W(K_x)-W(K) = \frac{2}{\omega(S^{n-1})}\ \int_{K^\circ
\setminus K_x^\circ } \|\xi\|^{-(n+1)} d\xi.
\end{equation}
The  bodies $K[t]$  have been
used by several authors (e.g. by B\"or\"oczky and Schneider
\cite{BoSch} and Glasauer and Gruber \cite{GG}) in connection with approximation of convex bodies 
by polytopes. We  generalize them as follows.
\vskip 3mm 
\noindent 
Let $f: K^\circ \rightarrow \mathbb{R}$ be a positive,  integrable function.
As  above, with $K_x$ instead of $M$,  we put
\begin{equation}\label{f-width}
w_f(x)=W_f(K_x)-W_f(K)= \frac{2}{\omega(S^{n-1})}\ \int_{K^\circ \setminus
K_x^\circ } f(\xi) d \xi
\end{equation}
and generalize (\ref{kt}) to
\begin{equation}\label{fkt}
K_f[t]=\{ x\in \mathbb{R}^n: w_f(x) \leq t\}.
\end{equation}
\par
\noindent
Thus, for instance, for $\beta \in \mathbb{R}$ and $f_\beta(\xi) = \|\xi\|^{-\beta}$ we get 
\begin{equation}\label{betakt}
K_{f_\beta}[t]=\left\{ x\in \mathbb{R}^n: \frac{2}{\omega(S^{n-1})}\ \int_{K^\circ \setminus
K_x^\circ } \|\xi\|^{-\beta} dx \leq t \right\},
\end{equation}
which,  in the particular case $\beta = n+1$,  gives the bodies (\ref{kt})  above.
\vskip 3mm
As $K_x =[x,K]$, $K_x^\circ= K^\circ \cap \{y \in \mathbb{R}^n: \langle y,x \rangle  \leq 1\}$.
Thus, putting $H^+\left(\frac{x}{\|x\|^2}, x\right)=  \{y \in \mathbb{R}^n: \langle y,x \rangle  \leq 1\}$, $K_x^\circ$ is obtained from $K^\circ$ by cutting off a cap 
$K^\circ \cap H^-\left(\frac{x}{\|x\|^2}, x\right)$ of $K^\circ$: 
$$
 K_x^\circ = K^\circ \cap H^+\left(\frac{x}{\|x\|^2}, \frac{x}{\|x\|}\right).
$$
and
$$
K^\circ\setminus K_x^\circ = K^\circ \cap H^-\left(\frac{x}{\|x\|^2}, \frac{x}{\|x\|}\right).
$$
Therefore
\begin{eqnarray} \label{form1}
K_{f}[t] &=& \left\{x \in \mathbb{R}^n: \frac{2}{\omega(S^{n-1})} \int_{ K^\circ \cap H^-\left(\frac{x}{\|x\|^2}, \frac{x}{\|x\|}\right)} f(\xi) d\xi \leq t \right\}.
\end{eqnarray}

\vskip 4mm 
\noindent
{\bf Remarks 1: 
 Properties of $K_f[t]$}
\vskip 2mm
\noindent
(i) It is clear that  for all $f$ and for all $t \geq 0$, $K \subset K_f[t]$ and that $K_{f_\beta}[0]=K$ for all $\beta$.
However, it can happen that $K$ is a proper subset of $K_f[0]$. 
\vskip 2mm
To see that, let $K=B^n_\infty =\{(x_1, \dots , x_n) \in \mathbb{R}^n: \mbox{max}_{1 \leq i \leq n} |x_i| \leq 1 \}$. Then $K^\circ =B^n_1=\{(x_1, \dots, x_n) \in \mathbb{R}^n: \sum_{i=1}^n
 |x_i| \leq 1 \}$.\\
Define $f:  B^n_1 \rightarrow \mathbb{R}$,  $(x_1, \dots, x_n) \rightarrow f( (x_1, \dots, x_n))$ by
\begin{equation*}
f(x)=\left\{
\begin{array}{cc}
0, &  x_n \geq 0\\
1, & ~\mbox{otherwise}.~
\end{array} \right.
\end{equation*}
Then $ (0, \dots, 0,\frac{3}{2}) \in K_f[0]$ but $ (0, \dots, 0, \frac{3}{2}) \notin K$.
\vskip 3mm
\noindent
(ii) $K_f[t]$ need neither  be bounded nor convex. Indeed, let  $K=B^2_\infty$.
Define $f:  B^2_1 \rightarrow \mathbb{R}$,  $(x_1, x_2) \rightarrow f( (x_1, x_2))$ by
\begin{equation*}
f(x)=\left\{
\begin{array}{cc}
\frac{1}{2}, &  x_2 \geq 0\\
1, & ~\mbox{otherwise}.~
\end{array} \right.
\end{equation*}
If $t \geq \frac{1}{\pi}$, $K_f[t]= \mathbb{R}^2$.  If $\frac{3}{4 \pi } \leq t <\frac{1}{\pi}$, $\{(x_1, x_2) \in \mathbb{R}^2: x_2 \geq 0\} \subset K_f[t]$. If $\frac{1}{2\pi } \leq t < \frac{3}{4\pi}$, $\{(0, x_2) \in \mathbb{R}^2: x_2 \geq 0\} \subset K_f[t]$. Thus $K_f[t]$ is unbounded in those cases. If $t < \frac{1}{2\pi}$, then $K_f[t]$ is bounded.
\par
Moreover, with the same $K$ and $f$: $\{(x_1, x_2) \in \mathbb{R}^2: x_2 \geq 0\} \subset K_f[ \frac{3}{4\pi}]$ and $\left(0, - \frac{1}{1-\sqrt{3}/2}\right) \in K_f[ \frac{3}{4\pi}]$. Let $x_0= \left( \frac{1}{1-\sqrt{3}/2}, \frac{-1}{1-\sqrt{3}/2} \right)$. Then $w_f(x_0) = \sqrt{3}\left(1-\sqrt{3}/16\right) > \frac{3}{4 \pi}$. Therefore, $K_f[ \frac{3}{4\pi}]$ is not convex.

\vskip 3mm
\noindent
(iii)  
Formulas (\ref{f-width}) and (\ref{form1})  show
that to define $K_f[t]$,  we cut off a set of  ``weighted
volume" $t$ of $K^\circ$. 
Thus $K_f[t]$  resembles  the convex floating body of $K^\circ$.
\par
Recall that for  $0 \leq \delta  \leq \frac{|K|}{2}$, the {\em convex floating body} $K_\delta$ of $K$ is the 
intersection of all halfspaces $H^+$ whose
defining hyperplanes $H$ cut off a set of volume at most $\delta$
from $K$ \cite{SW1}:
$$
K_{\delta}=\bigcap_{|H^-\cap K|\leq \delta} {H^+}. 
$$
For $\beta=0$, we get in formula (\ref{form1}), 
\begin{eqnarray*}
K_{f_0}[t] &=& \{x \in \mathbb{R}^n: \frac{2}{\omega(S^{n-1})} \int_{K^\circ \cap H^-\left(\frac{x}{\|x\|^2}, \frac{x}{\|x\|}\right)} d\xi \leq t \}\\
&=& \left\{x \in \mathbb{R}^n: \left|K^\circ \cap H^-\left(\frac{x}{\|x\|^2}, \frac{x}{\|x\|}\right) \right| \leq \frac{t \omega(S^{n-1} )}{2}\right\}
\end{eqnarray*}
\vskip 2mm
\noindent
However, $K_{f_0}[t]$ is not a convex floating body of $K^\circ$. 
\par
Indeed, it is easy to see that for the Euclidean ball $B=r B^n_2$ in $\mathbb{R}^n$ with radius $r$, 
 $B_{f_0}[t]$, for small $t$, is a Euclidean ball with radius of order
 $$
 r\left(1+ k_n r^\frac{2n}{n+1} t^\frac{2}{n+1}\right),
 $$ 
 where $k_n=\frac{1}{2} \left( \frac{n(n+1) |B^n_2|}{2 |B^{n-1}_2|}\right)^\frac{2}{n+1}$.
$(B^\circ)_\delta$, for small $\delta$,  is a ball with radius of order
$$
\frac{1}{r}\left(1-c_n r^\frac{2n}{n+1}  \delta^\frac{2}{n+1}\right),
$$ where $c_n=\frac{1}{2} \left( \frac{n+1}{ |B^{n-1}_2|}\right)^\frac{2}{n+1}$ (see e.g. \cite{SW1}) and 
$B_\delta$, for small $\delta$,  is a ball with radius of order
$$
r \left(1-\frac{c_n}{ r^\frac{2n}{n+1}}  \delta^\frac{2}{n+1}\right),
$$ (see also e.g. \cite{SW1}).
\vskip 2mm
Also, $K_{f_0}[t] $ is different from the {\em illumination body} $K^\delta$ which, for $\delta \geq 0$,   is defined as follows \cite{W1}:
$$
K^\delta= \{ x \in { \bf R}^n: |\mbox{co} [x,K] \backslash K| \leq \delta \}.
$$
Again, this can be seen by considering the  Euclidean ball $r B^n_2$. $(r B^n_2)^\delta$, for small $\delta$, is a Euclidean ball  
with radius of order
$$
r\left(1+ \frac{d_n}{ r^\frac{2n}{n+1} } \delta^\frac{2}{n+1}\right),
$$ where $d_n=\frac{1}{2} \left( \frac{n(n+1)}{ |B^{n-1}_2|}\right)^\frac{2}{n+1}$ \cite{W1}.

\vskip 4mm 

We have seen that $K_f[t]$ need not be convex. But it is always star-convex.   
\vskip 3mm 
\noindent
\begin{lemma} \label{starconvex}
Let $K$ be a convex body in $\mathbb{R}^n$ such that
$0$ is the center of gravity of $K$. Let $f: K^\circ \rightarrow \mathbb{R}$ be a positive, integrable function. 
\vskip 2mm 
\noindent
(i)  $K_f[t]$ is star convex i.e. $[0, x] \subset K_f[t]$ for all $x \in K_f[t]$.
\vskip 2mm 
\noindent
(ii) $K_f[t] = \bigcap _{s>0} K_f[t+s]$.
\end{lemma}
\vskip 3mm
\noindent
{\bf Proof.}
(i) Let $x \in K_f[t]$ and let $y \in [0,x]$.
Then $K_y =[y,K] \subset [x,K] = K_x$ and consequently $K^\circ \setminus K_y^\circ \subset K^\circ \setminus K_x^\circ$. As $f \geq 0$ on $K^\circ$, we therefore get
$$
\frac{2}{\omega(S^{n-1})} \int_{K^\circ\setminus K_y^\circ} f(\xi) d\xi \leq 
\frac{2}{\omega(S^{n-1})} \int_{K^\circ\setminus K_x^\circ} f(\xi) d\xi \leq t
$$
and thus $y \in K_f[t]$.
\vskip 2mm 
\noindent
(ii) For all $s>0$, $K_f[t] \subset K_f[t+s]$. Therefore, we only need to show that $\bigcap _{s>0} K_f[t+s] \subset K_f[t]$. Let thus  $x \in \bigcap _{s>0} K_f[t+s]$. Then for all $s>0$, $w_f(x) \leq t+s$. Letting $s \rightarrow 0$, we get $w_f(x) \leq t$. 
\vskip 4mm 
Additional
conditions on $f$ ensure convexity of $K_f[t]$. This is shown in the next lemma whose proof is the same as the corresponding 
one in \cite{BoSch}.
\vskip 3mm 
\noindent
\begin{lemma} \label{lemma:convex}
Let $K$ be a convex body in $\mathbb{R}^n$ such that
$0$ is the center of gravity of $K$. Let $f: S^{n-1} \rightarrow \mathbb{R}$ be a positive, integrable 
 function that is homogeneous of degree $\alpha$. Then
 $K_f[t]$ is convex for all  $\alpha \leq  -(n+1)$.

\end{lemma}
\vskip 3mm
\noindent
{\bf Proof.}
Let $x$ and $y$ be in $K_f[t]$ and let $0 <\lambda
<1$. For $t \in \mathbb{R}$, $t \geq 0$, the function $g(t)=t^{\gamma}$ is convex if $\gamma \geq 1$.
Therefore, and as $K_{(1-\lambda)x + \lambda y}
\subseteq (1-\lambda) K_x + \lambda K_y$, we get for $\alpha \leq -(n+1)$
$$
\frac{ h_{K_{(1-\lambda)x + \lambda y}}^{-(\alpha+n)}}{-(\alpha+n)}  \leq \frac{\left((1-\lambda)\
h_{K_{x}} + \lambda\  h_{K_{y}}\right)^{-(\alpha+n)} }{-(\alpha+n)} \leq \frac{(1-\lambda)\
h_{K_{x}} ^{-(\alpha+n)} + \lambda\  h_{K_{y}}^{-(\alpha+n)}}
{-(\alpha+n)}.
$$
Hence for $\alpha \leq -( n+1)$,
\begin{eqnarray*}
&& \frac{2}{-(\alpha+n)}\  \int_{S^{n-1}} f(u) 
h_{K_{(1-\lambda)x + \lambda y}}^{-(\alpha+n)}(u) d\sigma (u)
 \\
&& \leq \frac{2}{-(\alpha+n)}\  \bigg[ (1-\lambda) \int_{S^{n-1}} f(u)
h_{K_{x}} ^{-(\alpha+n)} (u) d\sigma (u) \\
&& +  \lambda  \int_{S^{n-1}} f(u) h_{K_{y}}^{-(\alpha+n)} (u) d\sigma(u) \bigg]
 \\
&& \leq (1-\lambda)\bigg[ \frac{2}{-(\alpha+n)}\
\int_{S^{n-1}}  f(u) h_{K} ^{-(\alpha+n)}(u)  d\sigma(u)   +  t \bigg]  \\
&&\hskip 20mm +  \lambda
\bigg[\frac{2}{-(\alpha+n)}\   \int_{S^{n-1}} f(u) 
h_{K}^{-(\alpha+n)}(u)d\sigma (u) + t \bigg] \\
&&= \frac{2}{-(\alpha+n)}\  \int_{S^{n-1}}  f(u)  h_{K}
^{-(\alpha+n)} (u) d\sigma   +  t.
\end{eqnarray*}
\vskip 3mm
\noindent
{\bf Remark.} If $\alpha > -(n+1)$, then $K_f[t]$ need not be convex. An example is the cube in $\mathbb{R}^2$ and the $f$ given in Remark 1 (ii).

\vskip 4mm 

Now we give conditions that guarantee  that $K_f[t]$ is bounded.  
\vskip 3mm 
\noindent
\begin{lemma} \label{bounded}
Let $K$ be a convex body in $\mathbb{R}^n$ such that
$0$ is the center of gravity of $K$. Let $f: K^\circ \rightarrow \mathbb{R}$ be a strictly positive, integrable function. Then
\vskip 2mm 
\noindent
(i)  $K_f[0] =K$.
\vskip 2mm 
\noindent
(ii) There exists $t_0$ such that for all $t \leq t_0$, $K_f[t] $ is bounded.
\vskip 2mm 
\noindent
(iii) Let $t \leq t_0$, where $t_0$ is as in (ii). Then we have for all $x \in \partial K_f[t] $ that $w_f(x)=t$.

\end{lemma}

\vskip 3mm
\noindent
{\bf Proof.}
\par
\noindent
(i) We only have to show that $K_f[0] \subset K$. Let $x \in K_f[0]$. Then $w_f(x)=\frac{2}{\omega(S^{n-1})} \int_{K^\circ \setminus K_x^\circ} f(\xi) d\xi =0$. As $f >0$ on $K^\circ$, this can only happen if $m(K^\circ \setminus K_x^\circ)=0$. As $K_x^\circ \subset K^\circ$ is closed and convex,  this can only happen if 
$K_x^\circ=K^\circ$, or, equivalently, $K_x=K$, or $x \in K$. 
\vskip 2mm 
\noindent
(ii) This follows immediately from (i),  Lemma \ref{starconvex} (ii) and the fact that, as $K$ is a  convex body, there exists $\alpha >0$ such that 
\begin{equation}\label{alpha}
B^n_2(0, \alpha) \subset K \subset B^n_2\left(0, \frac{1}{\alpha}\right).
\end{equation}
As
$K =K_f[0] = \bigcap_{t>0} K_f[t]$,  there exists $t_0$ such that for all $t \leq t_0$, $K_f[t] \subset 2K \subset  B^n_2\left(0, \frac{2}{\alpha}\right)$.
\vskip 2mm 
\noindent
(iii) Let $t \leq t_0$ and let $x \in \partial K_f[t]$. Suppose $w_f(x) <t$. Let $ y \in \{ax: a \geq 1\}$. Then 
$K_x=[x,K] \subset K_y = [y,K]$, hence $K_y^\circ \subset K_x^\circ$ and therefore 
$\int_{K^\circ \setminus K_y^\circ} f(\xi) d\xi \geq \int_{K^\circ \setminus K_x^\circ} f(\xi) d\xi $.
As $f >0$ on $K^\circ$, we can choose $y=ax$ with $ a>1$ such that $\frac{2}{\omega(S^{n-1})} \int_{K^\circ \setminus K_y^\circ} f(\xi) d\xi = t$. This implies that $x \notin \partial K_f[t]$, a contradiction.

\vskip 5mm 
\section{Relative entropies of cone measures and affine surface areas}
In this section we present new geometric interpretations of  important affine invariants, namely the $L_p$-affine surface areas.
Many such geometric interpretations have been given (see e.g. \cite{MW2, SW4, SW2004, W2, WY, WY1}). The remarkable fact here is that these geometric interpretations of affine invariants for {\em convex} bodies  are expressed in terms of {\em not necessarily convex} bodies,
a phenomenon which already occurred in \cite{WY1}.
\par
We also give new geometric interpretations for 
the relative entropies of cone measures of convex bodies. Geometric interpretations for those quantities were given first
in \cite{PW1}  in terms of $L_p$-centroid bodies:
For a convex body  $K$  in $\mathbb R^n$ of volume $1$ and $1\leq p \leq \infty$,  the $L_{p}$-centroid body $Z_{p}(K)$ is  this convex body that has support function 
\begin{equation*} \label{def:Zp}
h_{Z_{p}(K)}(\theta) = \left(\int_{K}|\langle x, \theta\rangle |^{p} dx \right)^{1/p}.
\end{equation*}
However, in the context of the  $L_p$-centroid bodies, the relative entropies appeared only after performing  a second order expansion of certain expressions. Now, using the mean width bodies, already a first order expansion makes them appear. Thus,  these bodies detect ``faster" more detail  of the  boundary of a convex body than the  $L_p$-centroid bodies.

\vskip 3mm \noindent 
\begin{theorem} \label{theorem1} 
Let $K$ be a convex body in $\mathbb{R}^n$ that is in $C^2_+$ and such that
$0$ is the center of gravity of $K$. Let $f:  K^\circ \rightarrow \mathbb{R}$ 
 be a continuous function such that
$f(y)\geq c$ for all $y\in  K^\circ$ and some constant $c>0$.
Then 
 \begin{equation}\nonumber
 \lim
_{t\rightarrow 0}  \frac{|K_f[t]|-|K|}{k_n \   t^{\frac{2}{n+1}}}=
 \int_{\partial K} \frac{\langle x, N_{K}(x) \rangle ^2 d \mu_K(x)}{f(y(x)) \kappa_K(x)^\frac{1}{n+1} }.
\end{equation}
$k_n=\frac{1}{2} \left( \frac{n(n+1) |B^n_2|}{2 |B^{n-1}_2|}\right)^\frac{2}{n+1}$  and $y(x) \in \partial K^\circ$ is such that $\langle y(x), x \rangle =1$.
 \end{theorem}
\vskip 3mm
\noindent {\bf Remark.} 
\par
\noindent
We put  $N_K(x)=u$. Then  $\langle x, N_{K}(x)\rangle = h_K(u)$ and $y(x)= \frac{u}{h_K(u)}$. As   $d\mu_K= f_K d \omega$, we therefore also have
 \begin{equation}\label{theorem1,sphere}
 \lim
_{t\rightarrow 0}  \frac{|K_f[t]|-|K|}{k_n \   t^{\frac{2}{n+1}}}=\int_{S^{n-1}}
\frac{ h_K(u)^{2} d\omega
(u) } {f_K(u)^\frac{n+2}{n+1} 
f \left(\frac{u}{h_K(u)}\right) }.
\end{equation}

\vskip 4mm
Theorem \ref{theorem1} leads to the announced new geometric interpretations of the above mentioned quantities
which  we introduce now.
\vskip 3mm
$L_p$-affine surface area, an extension of affine surface area, 
was introduced by Lutwak in the ground
breaking paper \cite{Lu2}  for $p >1$ and for general $p$ by Sch\"utt and Werner \cite{SW2004}.
For real  $p \neq -n$, we define  the
$L_p$-affine surface area $as_{p}(K)$ of $K$ as in \cite{Lu2} ($p
>1$) and \cite{SW2004} ($p <1, p \neq -n$) by
\begin{equation} \label{def:paffine}
as_{p}(K)=\int_{\partial K}\frac{\kappa_K(x)^{\frac{p}{n+p}}}
{\langle x,N_{ K}(x)\rangle ^{\frac{n(p-1)}{n+p}}} d\mu_{ K}(x) 
\end{equation}
and
\begin{equation}\label{def:infty}
as_{\pm\infty}(K)=\int_{\partial K}\frac{\kappa _K (x)}{\langle
x,N_{K} (x)\rangle ^{n}} d\mu_{K}(x), 
\end{equation}
provided the above integrals exist.
In particular, for $p=0$
$$
as_{0}(K)=\int_{\partial K} \langle x,N_{ K}(x)\rangle
\,d\mu_{K}(x) = n|K|.
$$
The case $p=1$ is the classical affine surface area which is independent
of the position of $K$ in space and which goes
back to Blaschke. 
$$
as_{1}(K)=\int_{\partial K} \kappa_{ K}(x)^\frac{1}{n+1} 
\,d\mu_{K}(x).
$$

Originally a basic affine invariant
from the field of affine differential geometry, it has recently
attracted increased attention too (e.g. 
\cite{LR1,  Lu2, MW1, SW1,  W1}).
\vskip 4mm
Then we have
\vskip 3mm 
\noindent 
\begin{corollary} \label{cor:p} 
Let $K$ be a convex body in $\mathbb{R}^n$ that is in $C^2_+$ and such that
$0$ is the center of gravity of $K$. 
\vskip 2mm
\noindent
(i) For $p \in \mathbb{R}$, $p \neq -n$, let $p_{as}: \partial K^\circ \rightarrow \mathbb{R}$ be defined by
$$
p_{as}(y) =\left( \frac{\langle x, N_{K}(x) \rangle }{\kappa_K(x)^\frac{1}{n+1}}\right)^\frac{n+p(n+2)}{n+p},
$$
where, for $y \in \partial K^\circ$, $x=x(y) \in \partial K$ is such that $\langle x, y \rangle=1$
Then 
 \begin{equation*}\label{}
  \lim
_{t\rightarrow 0}  \frac{|K_{p_{as}}[t]|-|K|}{
k_n \ t^{\frac{2}{n+1}}}=\int_{\partial K}
\frac{\kappa_K(x)^\frac{p}{n+p} d \mu_K(x)} {\langle x, N_{K}(x) \rangle ^\frac{n(p-1)}{n+p}}
= as_{p} (K).
\end{equation*}
\vskip 2mm
\noindent
(ii) For $\beta \in \mathbb{R}$, let $f_\beta: K^\circ \rightarrow \mathbb{R}$ be defined by
$$
f_\beta(y) =\frac{1}{\|y\|^\beta} = \langle x, N_{K}(x) \rangle ^\beta,
$$
where, again, 
for $y \in \partial K^\circ$, $x=x(y) \in \partial K$ is such that $\langle x, y \rangle=1$
Then 
 \begin{equation*}\label{}
  \lim
_{t\rightarrow 0}  \frac{|K_{f_\beta}[t]|-|K|}{
k_n\  t^{\frac{2}{n+1}}}=\int_{\partial K}
\frac{d \mu_K(x)}{\kappa_K(x)^\frac{1}{n+1} \langle x, N_{K}(x) \rangle ^{\beta-2}}
\end{equation*}

\end{corollary}

\vskip 2mm
\noindent
{\bf Proof.}
As $\partial K$ is in  $C^2_+$, the functions $p_{as}$ and $f_\beta$ 
satisfy the conditions of Theorem \ref{theorem1}. The proof  of the corollary  then follows immediately from Theorem \ref{theorem1}.
\vskip 4mm
\noindent 
{\bf Remarks} 
\vskip 2mm
\noindent
(i) For $\beta=0$,   we get in Corollary \ref{cor:p}  (ii) the $as_{-\frac{n}{n+2}}$-affine surface area of $K$.
\vskip 2mm
\noindent
(ii)  As $\kappa_K (r x) = r ^{-(n-1)} \kappa_K(x)$, it makes  most sense to put $f_{K}(ru) = f_{rK}(u)= r^{n-1}  f_K(u)$ and define $n-1$ to be  the degree of homogeneity of the function $f_K$.
Then $p_{as}$ is homogeneous of degree $\frac{2n(n+p(n+2))}{(n+1)(n+p)}$ and 
$f_\beta$ is homogeneous of degree $\beta$.
Thus, by Lemma \ref{lemma:convex}, 
$K_{p_{as}}[t]$ is convex if $-n < p \leq -n \frac{(n+1)^2+1}{(n+1)^2 +n+2}$ and $K_{f_\beta}[t]$ is convex if $\beta \leq -(n+1)$.

\vskip 4mm
Let $K$ a  convex body in $\mathbb R^n$ that is $C^{2}_{+}$. 
Let 
\begin{equation}\label{densities}
p_K(x)= \frac{ \kappa_{K}(x)}{\langle x, N_{K}(x) \rangle^{n} \  n|K^{\circ}|} \, , \   \ q_K(x)= \frac{\langle x, N_{K}(x) \rangle }{n\  |K|}.
\end{equation}
\noindent 
Then 
\begin{equation}\label{PQ}
P_K=p_K\  \mu_K \ \ \ \text{and}   \ \ \   Q_K=q_K \ \mu_K
\end{equation}
 are probability measures on $\partial K$ that are absolutely continuous with respect  to $\mu_{K}$.
\par
Recall now that the normalized cone measure $cm_K$
on $\partial K$ is defined as follows:
For every measurable set $A \subseteq \partial K$
\begin{equation}\label{def:conemeas} 
cm_{K}(A)  = \frac{1}{|K|}|\{ta : \ a \in A, t\in [0,1] \}|.
\end{equation}
\vskip 3mm
\noindent
The next proposition is well known. See e.g. \cite {PW1} for a proof. It shows that the measures $P_K$ and $Q_K$ defined in (\ref{PQ})
 are the cone measures  of $K^\circ$  and $K$.  $N_K:\partial K \rightarrow S^{n-1}$, $x \rightarrow N_K(x)$  is the Gauss map.
\vskip 3mm
\noindent 
\begin{proposition} \label{prop:conemeas}
\noindent 
Let $K$ a  convex body in $\mathbb R^n$ that is $C^{2}_{+}$. Let $P_K$ and $Q_K$ be the probability measures on $\partial K$  defined by (\ref{PQ}).  Then
$$
P_K= N_{K}^{-1}N_{K^{\circ}}cm_{K^{\circ}}\  \   \mbox{and} \   \ Q_K= cm_{K},
$$
\noindent 
or, equivalently, for every measurable subset $A$ in $ \partial K$
$$ 
P_K(A)= cm_{ K^{\circ}} \bigg(N_{{K^{\circ}}}^{-1} \big(N_{K} (A)\big)\bigg) \   \ \mbox{and} \   \   Q_K(A)= cm_{ K}(A).
$$
\end{proposition}
\vskip 4mm
In  the next two corollaries  we also use the following notations. For a convex body $K$ in $\mathbb{R}^n$ and $x \in \partial K$,
 let $r_i(x)$, $1 \leq i \leq n-1$ be the principal radii of curvature. We put
\begin{equation} \label{rR}
r= \text{inf}_{x \in \partial K} \min_{1 \leq i \leq n-1}  r_i(x) 
\ \  \text{and}  \ \  R= \sup_{x \in \partial K} \max_{1 \leq i \leq n-1}  r_i(x).
\end{equation}    
Note that if  $K$ be a convex body in $\mathbb{R}^n$ that is in $C^2_+$, then  $0 < r \leq R < \infty$.
Note also that $r=R$ iff $K$ is a Euclidean ball with radius $r$.
\vskip 4mm

\begin{corollary} \label{cor:omega} 
Let $K$ be a convex body in $\mathbb{R}^n$ that is in $C^2_+$ and such that
$0$ is the center of gravity of $K$.  Let $r, R$ be as in (\ref{rR}).
\vskip 2mm
\noindent
(i) 
Let $ent_1: \partial K^\circ \rightarrow \mathbb{R}$ be defined by
$$
ent_1(y) = 
\frac
{\kappa_K(x)^{-\frac{n+2}{n+1}}  \langle x, N_{K}(x) \rangle^{n+1}  }
{ \log \left(
\frac{ R^{2n} |K| \ \kappa_K(x)} {r^{2n} |K^\circ| \  \langle x, N_{K}(x) \rangle^{n+1}} \right) },
$$
where, again, for $y \in \partial K^\circ$, $x=x(y) \in \partial K$ is such that $\langle x, y \rangle=1$
Then 
 \begin{eqnarray*}\label{}
 \lim
_{t\rightarrow 0}  \frac{|K_{ent_1}[t]|-|K|}{
k_n\  t^{\frac{2}{n+1}}}&=& \int_{\partial K}
\frac{\kappa_K(x)} { \langle x, N_{K}(x) \rangle ^n} \log \frac{ R^{2n}  |K|\kappa_K(x)}{ r^{2n}   |K^\circ|\langle x, N_{K}(x) \rangle ^{n+1}} d \mu_K(x)\\
&=& n |K^\circ| \left[  [D_{KL}(P_K\| Q_K)  + 2n \log \left(\frac{R}{r} \right) \right]\\
&=& n |K^\circ|  \left[D_{KL}\big(N_{K}N_{K^{\circ}}^{-1}cm_{ K^{\circ}}\| cm_{ K}\big) + 2n \log \left(\frac{R}{r} \right) \right].
\end{eqnarray*}
\vskip 2mm
\noindent
\noindent
(ii)  Let $ent_2: \partial K^\circ \rightarrow \mathbb{R}$ be defined by
$$
ent_2(y) = 
\frac
{\kappa_K(x)^{-\frac{1}{n+1}}   }
{ \log \left(
\frac{R^{2n}  |K| \kappa_K(x)} { r^{2n} |K^\circ| \langle x, N_{K}(x) \rangle^{n+1}} \right) },
$$
where, again, for $y \in \partial K^\circ$, $x=x(y) \in \partial K$ is such that $\langle x, y \rangle=1$
Then 
 \begin{eqnarray*}\label{}
  \lim
_{t\rightarrow 0}  \frac{|K_{ent_2}[t]|-|K|}{
k_n\  t^{\frac{2}{n+1}}}&=&-  \int_{\partial K} \langle x, N_{K}(x) \rangle \log \frac{  r^{2n}  |K^\circ|\langle x, N_{K}(x) \rangle ^{n+1}}{  R^{2n}  |K| \kappa_K(x)}
 d \mu_K(x) \\
&= &- n |K| \left[ D_{KL}(Q_K || P_K) - 2n \log \left(\frac{R}{r} \right) \right]\\
&=& - n |K| \left[D_{KL}\big(cm_{ K} \| N_{K}N_{K^{\circ}}^{-1}cm_{ K^{\circ}}\big) - 2n \log \left(\frac{R}{r} \right)\right].
\end{eqnarray*}

\end{corollary}

\vskip 2mm
\noindent
{\bf Proof.}
As $\partial K$ is in  $C^2_+$,  $0 < r \leq R < \infty$ and  we have for all $x \in \partial K$ that
\begin{equation*}\label{condition}
B^n_2(x-r N_K(x), r) \subset K \subset  B^n_2(x-R N_K(x), R).
\end{equation*}
Suppose first that $r=R$. Then $K$ is a Euclidean ball with radius $r$ and  the right hand sides of the identities in the corollary are equal to $0$. Moreover, in this case, $ent_1$ and   $ent_2$  are identically equal to $\infty$. Therefore, for all $t\geq 0$, $K_{ent_1}[t]=K$
and $K_{ent_2}[t]=K$ and hence for all $t\geq 0$, $|K_{ent_1}[t]|-|K| =0$ and $|K_{ent_2}[t]|-|K| =0$. Therefore, the corollary holds trivially in this case.
\par
Suppose now that $r <R$. Then, as
$$
1 \leq \frac{ R^{2n}|K| \ \kappa_K(x)} {r^{2n} |K^\circ| \  \langle x, N_{K}(x) \rangle^{n+1}}   \leq\left( \frac{R}{r}\right)^{4n}.
$$
we get for all $x \in \partial K$ that
$$
f_{PQ}(x) \geq 
\left(\frac { |K^\circ|  r^{n-1}}{2  \log\left(\frac{R}{r}\right)}\right)^\frac{n-1}{2} >0.
$$
Thus  the functions $ent_1$ and   $ent_2$
satisfy the conditions of Theorem \ref{theorem1}. The proof  of the corollary  then follows immediately from Theorem \ref{theorem1}.

\vskip 4mm
In \cite{PW1}, the following   new affine invariant $\Omega_K$  was introduced and its relation to the relative entropies was established.
\vskip 2mm
\noindent
\par
\noindent
Let $K$ a convex body in $\mathbb R^n$ with centroid at the origin. 
$$ 
\Omega_{K} = \lim_{p\rightarrow \infty} \left(\frac{as_{p}(K)}{n |K^{\circ}|}\right)^{n+p}.
$$

\vskip 3mm
Let $p_K$ and $q_K$ be the densities defined in (\ref{densities}).
It was proved in \cite{PW1} that 
for  a  convex body $K$ in $\mathbb R^n$ that is $C^{2}_{+}$.
\begin{equation}\label{prop3:eq1} 
D_{KL}(P_K\|Q_K) = \log{\left( \frac{|K|}{|K^{\circ}|} \Omega_{K}^{-\frac{1}{n}} \right) }
\end{equation}
and
\vskip 2mm
\noindent
\begin{equation}\label{prop3:eq2} 
D_{KL}(Q_K\|P_K) = \log{ \left(  \frac{|K^{\circ}|}{|K|} \Omega_{K^{\circ}}^{-\frac{1}{n}} \right) }.
\end{equation}
\noindent 

\vskip 4mm
In  \cite{PW1},  geometric interpretations in terms of $L_p$-centroid bodies were given in the case of symmetric convex bodies for the new 
affine invariants $\Omega_K$. These interpretations are in the spirit of Corollary \ref{cor:p}: As $p \rightarrow  \infty$, the quantities  $\Omega_K$ and the related relative entropies appear in  appropriately chosen volume differences of  $K$ and its $L_p$-centroid bodies. However, in the context of the  $L_p$-centroid bodies,   a second order expansion was needed for the volume differences in order to make these terms appear.
Now,  it follows from Corollary \ref{cor:omega} (i) and  (ii)  
and Corollary \ref{cor:2} that no symmetry assumptions are needed and 
that  already a first order expansion gives such  geometric interpretations,  if one uses the mean width bodies instead of the 
$L_p$-centroid body.
\vskip 4mm 
\noindent 
\begin{corollary} \label{cor:2} 
Let $K$ be a convex body in $\mathbb{R}^n$ that is in $C^2_+$ and such that
$0$ is the center of gravity of $K$. Let the functions $ent_1$ and $ent_2$ be as in Corollary \ref{cor:p}.
Then
 \begin{equation*}\label{}
  \lim
_{t\rightarrow 0}  \frac{|K_{ent_1}[t]|-|K|}{
k_n\  t^{\frac{2}{n+1}}} - 2 n^2 |K^\circ|  \log \left(\frac{R}{r} \right) = n |K^\circ| \log \left( \frac{|K|}{|K^\circ|} \Omega_{K}^{-\frac{1}{n}} \right).
\end{equation*}
and 
 \begin{equation*}\label{}
  \lim
_{t\rightarrow 0}  \frac{|K_{ent_2}[t]|-|K|}{
k_n\  t^{\frac{2}{n+1}}} - 2 n^2|K|  \log \left(\frac{R}{r} \right)= n |K| \log \left( \frac{|K|}{|K^\circ|} \Omega_{K^\circ}^\frac{1}{n} \right).
\end{equation*}

\end{corollary}
\vskip 4mm
\noindent

\section{Proof of Theorem \ref{theorem1}}

To prove Theorem \ref{theorem1}, we need the following lemmas. The first one, Lemma \ref{ellipse},  is well known.
\vskip 4mm
\noindent
\begin{lemma} \label{ellipse}
Let $\mathcal{E}_n(x_0,a)$ be an ellipsoid in $\mathbb{R}^n$ centered at $x_0$ and with axes parallel to the coordinate axes and of lengths $a_1, \dots, a_n$. Let $0 < \Delta < a_n$. Let 
$$
C(\mathcal{E}_n, \Delta) = \mathcal{E}_n \cap H(x_0+ (a_n-\Delta) e_n,e_n)
$$
be a cap of $\mathcal{E}_n(x_0,a)$ of height $\Delta$.
Then 
\begin{eqnarray*}
 \ \frac{2^\frac{n+1}{2}\left( 1- \frac{\Delta}{2a_n}\right)^\frac{n-1}{2}  \ |B^{n-1}_{2}|}{n+1}  \prod_{i=1}^{n-1} \frac{a_i}{\sqrt{a_n}} \  \Delta^\frac{n+1}{2}  
\leq
\left| C(\mathcal{E}_n, \Delta) \right| \\
\leq 
  \ \frac{2^\frac{n+1}{2} \ |B^{n-1}_{2}|}{n+1}  \prod_{i=1}^{n-1} \frac{a_i}{\sqrt{a_n}} \  \Delta^\frac{n+1}{2}  
\end{eqnarray*}
\end{lemma}

\vskip 4mm
In the next few lemmas and throughout the remainder of the paper  we will use the following notation.
\par
Let $K$ be a convex body in $\mathbb{R}^n$. Let  $f:  K^\circ \rightarrow \mathbb{R}$ 
 be an integrable function and for $t \geq 0$,  let $K_f[t]$ be a mean width body of $K$.
 For $x \in \partial K$, let
\begin{equation}\label{xt}
x_t = \{\gamma x: \gamma \geq 0 \} \cap \partial K_f[t].
\end{equation}
\vskip 2mm
Let $y(x) \in \partial K^\circ$ be such that $\langle y(x), x \rangle =1$.
Let $m$ be the Lebesgue measure on $\mathbb{R}^n$ and let $m_f$ be the measure (on $K^\circ$) defined by $m_f=\frac{2f}{\omega(S^{n-1})} \  m$, i.e.
for all $A \subset K^\circ$
\begin{equation*}
m_f(A)=\frac{2}{\omega(S^{n-1})} \int_Af(\xi) d \xi. 
\end{equation*}

\vskip 4mm
\noindent
\begin{lemma} \label{interchange}
Let $K$ be a convex body in $\mathbb{R}^n$ that is in $C^2_+$ and such that
$0$ is the center of gravity of $K$. Let $f:  K^\circ \rightarrow \mathbb{R}$ 
 be an integrable function such that
$f(y)\geq c$ for all $y\in  K^\circ$ and some constant $c>0$. Let $x_t$ be as in (\ref{xt}). Then
the functions 
$$
\frac{1}{t^\frac{2}{n+1}}\  \left(\frac{ \|x_t\| } {\|x\|} - 1\right)
$$
are uniformly (in $t$) bounded by an integrable  function.

\end{lemma}

\vskip 3mm
\noindent
{\bf Proof.}
We can assume that $t \leq t_0$ where $t_0$ is given by Lemma \ref{bounded}. Then $K_f[t]$ is bounded and hence 
\begin{equation}\label{enthalten}
K_f[t] \subset B^n_2(0,a)
\end{equation}
for some $a >0$. As $f \geq c$ on $K^\circ$, we get with (\ref{form1})
\begin{eqnarray*}
t &\geq&  \frac{2}{\omega(S^{n-1})} \int_{ K^\circ \cap H\left(\frac{x_t}{\|x_t\|^2}, \frac{x}{\|x\|}\right)^-} f(\xi) d\xi \\
&\geq& \frac{2 c}{\omega(S^{n-1})} \left|K^\circ \cap H^-\left(\frac{x_t}{\|x_t\|^2}, \frac{x}{\|x\|}\right)\right|.
\end{eqnarray*}
As $K$ is in $C^2_+$, $K^\circ$ is in $C^2_+$. Thus, by the Blaschke rolling theorem (see \cite{Schneider1993}), there exists $r_0 >0$ such that for all $y \in \partial K^\circ$, $B^n_2(y-r_0 N_{K^\circ}(y), r_0) \subset K^\circ$. Let now $y(x) \in \partial K^\circ$ be such that $\langle x,  y(x) \rangle=1$. Then $N_{K^\circ}(y(x))=\frac{x}{\|x\|}$ and thus 
\begin{eqnarray*}
t &\geq&  \frac{2 c }{\omega(S^{n-1})} \left|B^n_2\left(y(x)-r_0 \frac{x}{\|x\|}, r_0\right) \cap H^-\left(\frac{x_t}{\|x_t\|^2}, \frac{x}{\|x\|}\right)\right| \\
&\geq&  \frac{2^\frac{n+3}{2}\   c \  r_0^\frac{n-1}{2} \left|B^{n-1}_2\right| }{(n+1) \  \omega(S^{n-1})}  \left(\frac{1}{\|x\|} - \frac{1}{\|x_t\|}\right) ^\frac{n+1}{2}, 
\end{eqnarray*}
where we have used that  $ \left|B^n_2\left(y(x)-r_0 \frac{x}{\|x\|}, r_0\right) \cap H^-\left(\frac{x_t}{\|x_t\|^2}, \frac{x}{\|x\|}\right)\right|$ is the volume of a cap of height $ \frac{1}{\|x\|} - \frac{1}{\|x_t\|}=\frac{\|x_t - x\|}{\|x_t\| \|x\|}$ of the ball $B^n_2\left(y(x)-r_0 \frac{x}{\|x\|}, r_0\right)$ which we have  estimated from below using Lemma \ref{ellipse}. We assume also that $t$ is so small that  $ \frac{1}{\|x\|} - \frac{1}{\|x_t\|} < r_0$.
\par
As $x$ and $x_t$ are colinear, $\frac{ \|x_t\| } {\|x\|} - 1 = \frac{\|x_t - x\|}{\|x\|}$ and hence
\begin{eqnarray}\label{Bounded}
\frac{1}{t^\frac{2}{n+1}}\  \left(\frac{ \|x_t\| } {\|x\|} - 1\right) &=& \frac{1}{t^\frac{2}{n+1}}\   \frac{\|x_t - x\|}{\|x\|}
\leq \left(\frac{(n+1)\  \omega(S^{n-1})}{ c\  \left|B^{n-1}_2\right|}\right)^\frac{2}{n+1} 
\   \frac{r_0^{-\frac{n-1}{n+1}}}{2 ^\frac{n+3}{n+1}} \|x_t\| \nonumber \\
&\leq &  \left(\frac{(n+1)\  \omega(S^{n-1})}{ c\  \left|B^{n-1}_2\right|}\right)^\frac{2}{n+1} 
\   \frac{r_0^{-\frac{n-1}{n+1}}}{2^\frac{n+3}{n+1}}\  a.
\end{eqnarray}
In the last inequality we have used (\ref{enthalten}).  The expression (\ref{Bounded}) is a constant and thus integrable.

\vskip 4mm
\noindent
\begin{lemma} \label{limit}
Let $K$ be a convex body in $\mathbb{R}^n$ that is in $C^2_+$ and such that
$0$ is the center of gravity of $K$. Let $f:  K^\circ \rightarrow \mathbb{R}$ 
 be a continuous,  positive  function.
Then for  all $x \in \partial K$ one has 
\begin{eqnarray*}
\lim_{t \rightarrow 0} \frac{ \langle x, N_K(x) \rangle }{n \ k_n \  t^\frac{2}{n+1}}\   \left[\left(\frac{\|x_t\|}{\|x\|}\right)^n - 1  \right] =  \frac{\langle x, N_K(x) \rangle ^2}{\kappa_K(x)^\frac{1}{n+1} f(y(x))^\frac{2}{n+1}}, 
\end{eqnarray*}
where $k_n=\frac{1}{2} \left( \frac{n(n+1) |B^n_2|}{2 |B^{n-1}_2|}\right)^\frac{2}{n+1}$ and $y(x) \in \partial K^\circ$ is such that $\langle  x, y(x) \rangle=1$.
\end{lemma}

\vskip 2mm
\noindent
{\bf Proof.}
Let $x \in \partial K$. Let $x_t$ be as in (\ref{xt}).   
As $x$ and $x_{t}$ are collinear and as
$(1+s)^n\geq 1+n s$ for $s\in [0,1)$, one has for small enough $t$,
\begin{eqnarray*}
\frac{\langle x,
N_{K}(x)\rangle}{n}\left[\left(\frac{\|x_t\|}{\|x\|}\right)^n-1\right]=\frac{\langle
x,
N_{K}(x)\rangle}{n}\left[\left(1+\frac{\|x_t-x\|}{\|x\|}\right)^n-1\right]\geq
\Delta(x,t), 
\end{eqnarray*}
where $\Delta(x,t)=\left\langle \frac{x}{\|x\|}, N_{
K}(x)\right\rangle \|x_t-x\|=\langle x_t-x,N_K(x)\rangle$.
\par
\noindent
Similarly, as $(1+s)^n \leq 1+n s+2^n s^2$ for $s\in [0,1)$, one has for $t$ small enough,
\begin{equation}\label{Form:0:0}
\frac{\langle x,
N_{K}(x)\rangle}{n}\left[\left(\frac{\|x_t\|}{\|x\|}\right)^n-1\right]\leq
\Delta(x,t) \left[ 1+\frac{2^n}{n} \ \left(\frac{\|x_t
-x\|}{\|x\|}\right) \right].\
\end{equation}
\par\noindent
Hence for  $\varepsilon>0$ there exists $t_\varepsilon \leq t_0$, $t_0$ from Lemma \ref{bounded},  such that for all
$0<t\leq t_\varepsilon$
\begin{equation*}
1\leq \frac{\langle x,
N_{K}(x)\rangle\left[\left(\frac{\|x_t\|}{\|x\|}\right)^n-1\right]}{n\  \Delta(x,t)}\leq
1+\varepsilon.
\end{equation*}
By Lemma \ref{bounded} (iii), $m_f(K^\circ\setminus K_{x_t}^\circ) =t$ and thus 
\begin{equation*}
1\leq \frac{\langle x,
N_{K}(x)\rangle\left[\left(\frac{\|x_t\|}{\|x\|}\right)^n-1\right] \left(m_f(K^\circ\setminus K_{x_t}^\circ)\right)^\frac{2}{n+1}}{n\  \Delta(x,t)\   t^\frac{2}{n+1}}\leq
1+\varepsilon.
\end{equation*}
Let now 
 $y=y(x) \in \partial K^\circ$  be  such that $\langle x , y \rangle =1$. Thus 
$y = \frac{N_K(x)}{\langle x, N_K(x) \rangle}$ and $N_{K^\circ}(y)= \frac{x}{\|x\|}$.  
As $f$ is continuous on $K^\circ$, there exists $\delta >0$ such that for all $z \in B^n_2(y, \delta)$,
$$
f(y) -\varepsilon < f(z) < f(y) + \varepsilon.
$$
We choose $t$ so small that $K^\circ\setminus K_{x_t}^\circ \subset B^n_2(y, \delta)$.
Then 
\begin{eqnarray*}
&&\hskip - 25mm \frac{2 \left(f(y(x)) - \varepsilon\right) }{\omega(S^{n-1})} \left| K^\circ\setminus K_{x_t}^\circ\right| \leq \\
&&m_f\left (K^\circ\setminus K_{x_t}^\circ)\right) = \frac{2}{\omega(S^{n-1})} \int_{K^\circ \setminus K_{x_t}^\circ} f d\xi \\
&& \hskip 30mm \leq \frac{2 \left(f(y(x)) + \varepsilon\right) }{\omega(S^{n-1})} \left| K^\circ\setminus K_{x_t}^\circ\right|
\end{eqnarray*}
and  we get with (new) absolute constants $c_1$ and  $c_2$ that
\begin{eqnarray}\label{gleichung1}
1-c_1 \varepsilon &\leq& \frac{\langle x,
N_{K}(x)\rangle\left[\left(\frac{\|x_t\|}{\|x\|}\right)^n-1\right] \left(\frac{2 f(y(x))}{\omega(S^{n-1})} \left| K^\circ\setminus K_{x_t}^\circ\right|\right)^\frac{2}{n+1}}{n\  \Delta(x,t)\   t^\frac{2}{n+1}} \nonumber \\
&\leq&
1+c_2 \varepsilon.
\end{eqnarray}
\par
As $K$ and hence $K^\circ$  is in $C^2_+$, $ \kappa_{K^\circ}(y) >0$. 
It is well known  (see \cite{SW4}) that then there exists  an ellipsoid $\mathcal {E}=\mathcal {E}(y-a_n N_{K^\circ}(y), a)$ centered at $y-a_n N_{K^\circ}(y)$ and with half axes of lengths $a_1 \dots a_n$ which approximates $\partial  K^\circ$
in a neighborhood of $y$. For the computations that follow, we can assume without loss of generality that $N_{K^\circ}(y)=e_n$ and that the other axes of $\mathcal{E}$ coincide with $e_1 \dots, e_{n-1}$. Thus (see \cite{SW4}), 
for $\varepsilon >0$  given, 
there exists $\Delta_\varepsilon$ such that for all $\Delta \leq  \Delta_\varepsilon$ 
\begin{eqnarray} \label{E}
&&\hskip -20mm \mathcal {E}\big(y-(1-\varepsilon) a_n N_{K^\circ}(y), (1-\varepsilon)a \big)  \  \cap H_\Delta^-  \nonumber \\
&&\hskip 10mm \subseteq K^\circ \cap H_\Delta^- 
\subseteq \nonumber \\
&& \hskip 20mm \mathcal {E}\big(y-(1+\varepsilon)a_n N_{K^\circ}(y), (1+\varepsilon)a  \big) \cap H_\Delta^-,
\end{eqnarray}
where $H_\Delta=H(y-\Delta e_n,e_n)$.
Also (see  \cite{SW4}),  
\begin{equation}\label{kappa}
\kappa_{ K^\circ}(y)= \prod_{i=1}^{n-1} \frac{a_{n}}{a_i^2}.
\end{equation} 
As $x_t \rightarrow x$ as $t \rightarrow 0$, we can choose $t$ so small that 
$K^\circ\setminus K_{x_t}^\circ=K^\circ \cap H^-\left(\frac{x_t}{\|x_t\|^2}, \frac{x}{\|x\|}\right) $ is contained in 
$H^-(y-\Delta e_n,e_n)$.
Hence, by (\ref{E}), 
\begin{eqnarray*}
 \left|  \mathcal {E}\big(y-(1-\varepsilon) a_n N_{K^\circ}(y), (1-\varepsilon)a \big)  \cap H^-\left(\frac{x_t}{\|x_t\|^2}, \frac{x}{\|x\|}\right) \right|  \leq 
 \left| K^\circ\setminus K_{x_t}^\circ \right|  \leq \\
 \left|  \mathcal {E}\big(y-(1+\varepsilon) a_n N_{K^\circ}(y), (1+\varepsilon)a \big)  \cap H^-\left(\frac{x_t}{\|x_t\|^2}, \frac{x}{\|x\|}\right) \right|. 
 \end{eqnarray*}
By Lemma \ref{ellipse}, with  (\ref{kappa}),  and  as $\frac{1}{\|x\|} - \frac{1}{\|x_t\|}= \frac{\Delta(x,t)}{\|x_t\| \langle x, N_K(x) \rangle}$, we get  with new absolute constants $c_1$ and $c_2$
\begin{eqnarray*}
 (1- c_1 \varepsilon) 
 \frac{2^\frac{n+1}{2}   \left|B^{n-1}_2\right| }{(n+1)\left(\kappa_{K^\circ} (y)\right)^{\frac{1}{2}} }   
 \left(\frac{\Delta(x,t)}{\|x_t\| \langle x, N_K(x) \rangle}\right) ^\frac{n+1}{2}  
\leq  \left| K^\circ\setminus K_{x_t}^\circ \right|  \leq \\ (1+ c_2 \varepsilon)  \frac{2^\frac{n+1}{2}   \left|B^{n-1}_2\right| }{(n+1)\left(\kappa_{K^\circ} (y)\right)^{\frac{1}{2}} }    \left(\frac{1}{\|x\|} - \frac{1}{\|x_t\|}\right) ^\frac{n+1}{2} \\
  \hskip30mm = (1+ c_2 \varepsilon)  \frac{2^\frac{n+1}{2}   \left|B^{n-1}_2\right| }{(n+1)\left(\kappa_{K^\circ} (y)\right)^{\frac{1}{2}} } 
\left(\frac{\Delta(x,t)}{\|x_t\| \langle x, N_K(x) \rangle}\right) ^\frac{n+1}{2}.  
 \end{eqnarray*}
Hence, again with new absolute constants $c_1$ and $c_2$, (\ref{gleichung1}) becomes
\begin{equation*}\label{gleichung2}
1-c_1 \varepsilon \leq \frac{\langle x,
N_{K}(x)\rangle\left[\left(\frac{\|x_t\|}{\|x\|}\right)^n-1\right] 2\ \left(\frac{2  f(y) |B^{n-1}_2}{(n+1)\omega(S^{n-1})} \right)^\frac{2}{n+1}}{n\    t^\frac{2}{n+1} \left(\kappa_{K^\circ} (y)\right)^{\frac{1}{n+1}} \|x_t\|  \langle x, N_K(x) \rangle}\leq
1+c_2 \varepsilon.
\end{equation*}
Therefore, as $\|x_t\| \rightarrow \|x\|$ as $t \rightarrow 0$, 
\begin{eqnarray*}
&&\hskip -10mm \lim_{t \rightarrow 0} \frac{ \langle x, N_K(x) \rangle }{n  \  t^\frac{2}{n+1}}\   \left[\left(\frac{\|x_t\|}{\|x\|}\right)^n - 1  \right] =  \\ 
&&\hskip 20mm \frac{1}{2} \left( \frac{n(n+1) |B^n_2|}{2 |B^{n-1}_2|}\right)^\frac{2}{n+1}\   \frac{\kappa_{K^\circ}(y)^\frac{1}{n+1}  \|x\| \langle x, N_K(x) \rangle }{\ f(y)^\frac{2}{n+1}}
 \end{eqnarray*}
Now we use that $\|x\| =  \frac{1}{\langle y, N_{K^\circ} (y) \rangle }$ and that (see e.g. \cite{WY1})
$$
\frac{\kappa_{K^\circ}(y)^\frac{1}{n+1} }{\langle y, N_{K^\circ} (y) \rangle} = \frac{\langle x, N_{K} (x) \rangle}{\kappa_{K}(x)^\frac{1}{n+1} }
$$
We put  $k_n= \frac{1}{2} \left( \frac{n(n+1) |B^n_2|}{2 |B^{n-1}_2|}\right)^\frac{2}{n+1}$ and  get that 
$$
\lim_{t \rightarrow 0} \frac{ \langle x, N_K(x) \rangle }{n  \  t^\frac{2}{n+1}}\   \left[\left(\frac{\|x_t\|}{\|x\|}\right)^n - 1  \right] = k_n   \frac{ \langle x, N_K(x) \rangle ^2 }{\kappa_{K}(x)^\frac{1}{n+1} \ f(y)^\frac{2}{n+1}}.
$$

\vskip 5mm
\noindent
{\bf Proof of Theorem \ref{theorem1}}
\vskip 3mm
\noindent
It is  well known (see e.g. \cite{WY1}),  that
for a convex body  $K$  and a star-convex body $L$  
with  $0\in int(K)$ and $K\subset L$
$$
|L|-|K|=\frac{1}{n}\int _{\partial K} \langle x,
N_{K}(x)\rangle\left[\left(\frac{\|x'\|}{\|x\|}\right)^n-1\right]\,d\mu_{K}
(x)
$$
where $x\in \partial K$, $x'\in \partial L$ and  $ x=\partial
K\cap [0,x']$. \par
\noindent
Therefore, 
$$
|K_f[t]| -|K |= \frac{1}{n}\  \int_{\partial K } \langle x, N_K(x) \rangle \left(\left(\frac{\|x_t\|}{\|x\|}\right)^n - 1  \right) d  \mu_K(x)
$$
We now use Lemma \ref{interchange}  and Lebegue's theorem to interchange integration and limit and then Lemma \ref{limit} and get
\begin{eqnarray*} \label{diffvol}
 \lim _{t \rightarrow 0}  \frac{|K_f[t]| -|K |}{t^\frac{2}{n+1} }
 &= & \frac{1}{n} \lim _{t \rightarrow 0} \frac{1}{t^\frac{2}{n+1} } 
\  \int_{\partial K} \langle x, N_K(x) \rangle  \left[\left(\frac{\|x_t\|}{\|x\|}\right)^n - 1  \right] d  \mu_K(x) \\
&= & \int_{\partial K}\lim _{t \rightarrow 0}  \frac{\langle x, N_K(x) \rangle }{ n \  t^\frac{2}{n+1}}\   \left[\left(\frac{\|x_t\|}{\|x\|}\right)^n - 1  \right] d  \mu_K(x) \\
&= & k_n    \int_{\partial K} \frac{ \langle x, N_K(x) \rangle ^2 }{\kappa_{K}(x)^\frac{1}{n+1} \ f(y)^\frac{2}{n+1}} d  \mu_K(x).
\end{eqnarray*}
This finishes the proof of Theorem \ref{theorem1}.

 \vskip 2mm \noindent Elisabeth Werner\\
{\small Department of Mathematics \ \ \ \ \ \ \ \ \ \ \ \ \ \ \ \ \ \ \ Universit\'{e} de Lille 1}\\
{\small Case Western Reserve University \ \ \ \ \ \ \ \ \ \ \ \ \ UFR de Math\'{e}matique }\\
{\small Cleveland, Ohio 44106, U. S. A. \ \ \ \ \ \ \ \ \ \ \ \ \ \ \ 59655 Villeneuve d'Ascq, France}\\
{\small \tt elisabeth.werner@case.edu}\\ \\
\and Justin Jenkinson\\
{\small Department of Mathematics}\\
{\small Case Western Reserve University}\\
{\small Cleveland, Ohio 44106, U. S. A.}\\
{\small \tt jdj13@case.edu}


\newpage
 
 \begin{thebibliography}{~~}


\bibitem {BoSch}{\sc K. B\"or\"oczky Jr., R. Schneider},  {\em The mean width of circumscribed random polytopes}, Canadian Math. Bull. {\bf 53} (2010), 614-628. 


\bibitem {CT}{\sc T. Cover and J. Thomas}, {\em Elements of information theory}, second ed., Wiley-Interscience,  (John Wiley and Sons), Hoboken, NJ, (2006). 

\bibitem {DCT}
{\sc A. Dembo, T. Cover, and J. Thomas} {\em  Information theoretic inequalities},  IEEETrans. Inform. Theory {\bf 37} (1991), 1501-1518.







\bibitem{Ga1}
{\sc  R. J. Gardner}, { \em A positive answer to the Busemann-Petty problem in three dimensions}, 
Ann. of Math. {\bf 140} no.2 (1994), 435-47.  


\bibitem{Ga3}
{\sc R. J. Gardner}, {\em The dual Brunn-Minkowski theory for bounded Borel sets: Dual affine 
quermassintegrals and inequalities}, Adv. Math. {\bf 216} (2007), 358-386. 



\bibitem{GaKoSch}
{\sc  R. J. Gardner, A. Koldobsky, and T.  Schlumprecht}, {\em An analytical solution 
to the Busemann-Petty problem on sections of convex bodies}, Ann. of Math. {\bf 149} no.2 (1999), 691-703.



\bibitem{GaZ}
{\sc R. J. Gardner and G. Zhang},
{\em Affine inequalities and radial mean bodies.}
 Amer. J. Math. {\bf 120} no.3 (1998), 505-528.



\bibitem{GrZh}
{\sc
E. Grinberg and G. Zhang,} {\em Convolutions, transforms, and convex bodies}, Proc. 
London Math. Soc. {\bf 78} no.3 (1999), 77-115.

\bibitem{GG}
{\sc S. Glasauer and P. M. Gruber}, {\em Asymptotic estimates for best and stepwise approximation of
convex bodies III}, Forum Math. {\bf 9} (1997), 383-404.

\bibitem{Hab}
{\sc C. Haberl}, {\em Blaschke valuations}, Amer. J. of Math., in press



\bibitem{HabSch2}
{\sc C. Haberl and F.  Schuster,} {\em General Lp affine isoperimetric inequalities}.
J. Differential Geometry {\bf 83} (2009), 1-26.

\bibitem{HLYZ}
{\sc C. Haberl, E. Lutwak, D. Yang and G. Zhang,} {\em The even Orlicz Minkowski problem}, 
Adv. Math. {\bf 224} (2010), 2485-2510




\bibitem{Klain1}
{\sc D. Klain},
{\em  Star valuations and dual mixed volumes}, Adv. Math. {\bf 121} (1996), 80-101. 

\bibitem{Klain2}
{\sc D. Klain},
{\em Invariant valuations on star-shaped sets},   Adv. Math.  {\bf 125} (1997), 95-113. 




\bibitem{Lud2}
{\sc M. Ludwig}, {\em Ellipsoids and matrix valued valuations}, Duke Math. J. {\bf 119} (2003), 159-188.


\bibitem{Lud3}
{\sc M. Ludwig}, {\em
Minkowski areas and valuations}, 
J. Differential Geometry, {\bf 86} (2010), 133-162.

\bibitem{LR1}
{\sc M. Ludwig and M. Reitzner,}  {\em A Characterization of Affine Surface Area}, Adv.  Math. {\bf 147} (1999), 138-172.

\bibitem{LR2}
{\sc M. Ludwig and M. Reitzner,} {\em A classification of $SL(n)$
invariant valuations.}  Ann. of Math. {\bf 172 } (2010), 1223-1271. 



\bibitem{Lu1988}{\sc E. Lutwak}, {\em Intersection bodies and dual mixed volumes}, Adv.  Math. {\bf 71} (1988), 232-261.




\bibitem{Lu2}{\sc E. Lutwak}, {\em The Brunn-Minkowski-Firey theory II : Affine and
geominimal surface areas}, Adv. Math. {\bf 118}  (1996),   244-294.




\bibitem{LYZ2000} {\sc E. Lutwak, D. Yang and G. Zhang}, {\em A new ellipsoid associated with convex bodies}, Duke Math. J. { \bf 104} (2000), 375-390.

\bibitem{LZ} {\sc E. Lutwak and G. Zhang}, {\em Blaschke-Santal\'{o}
inequalities}, J. Differential Geom. {\bf 47} (1997), 1-16.


\bibitem{LYZ2002} {\sc E. Lutwak, D. Yang and G. Zhang}, {\em Sharp Affine $L_p$ Sobolev inequalities}, 
J. Differential Geometry {\bf 62} (2002), 17-38.


\bibitem{LYZ2002/1} {\sc E. Lutwak, D. Yang and G. Zhang}, {\em The Cramer--Rao inequality for star bodies}, Duke Math. J. {\bf 112} (2002), 59-81.

\bibitem{LYZ2004} {\sc E. Lutwak, D. Yang and G. Zhang}, {\em Volume inequalities for subspaces of $L_p$},  
J. Differential Geometry  {\bf 68} (2004), 159-184.

\bibitem{LYZ3} {\sc E. Lutwak, D. Yang and G. Zhang}, {\em Moment-entropy inequalities}, Ann. Probab. {\bf 32} (2004), 757Ð774. 

\bibitem{MW1}
{\sc M. Meyer and E. Werner}, {\em The
Santal${\mbox{\'o}}$-regions of a convex body.} Transactions of
the AMS {\bf 350} no.11, (1998) 4569-4591.


\bibitem{MW2}
{\sc M. Meyer and E. Werner}, {\em On the p-affine surface area.}
Adv. Math. {\bf 152} (2000),  288-313.

\bibitem{NPRZ}
{\sc F. Nazarov,  F. Petrov, D. Ryabogin and A. Zvavitch},  {\em A remark on the Mahler conjecture: local minimality of the unit cube}, 
Duke Math. J. {\bf 154}, (2010), 419-430.

\bibitem{PW1}
{\sc G. Paouris and E. Werner}, {\em Relative entropy of cone measures and $L_p$ centroid bodies}, 
preprint

\bibitem{RuZ}
{\sc B. Rubin and G. Zhang}, {\em Generalizations of the Busemann-Petty problem for sections of convex bodies}, J. Funct. Anal. { \bf 213} (2004), 473-501.





\bibitem{Schneider1993} 
{\sc R. Schneider},
{\em Convex Bodies: The Brunn-Minkowski
theory.} Cambridge University Press, (1993).

\bibitem{Schu}
{\sc
F. Schuster}, {\em Crofton measures and Minkowski valuations}, Duke Math. J. {\bf 154}, (2010), 1-30. 

\bibitem{SW1}
{\sc C. Sch\"{u}tt and  E. Werner}, {\em The convex floating
body.} Math. Scand. {\bf 66} (1990), 275-290.


\bibitem{SW4}
{\sc C. Sch{\"u}tt and E. Werner}, {\em Random polytopes of points
chosen from the boundary of a convex body.} GAFA Seminar Notes,
Lecture Notes in Mathematics, Springer-Verlag {\bf 1807} (2002), 241-422,
.



\bibitem{SW2004}
{\sc C. Sch{\"u}tt and E. Werner}, {\em Surface bodies and
p-affine surface area.} Adv. Math. {\bf 187}  (2004), 98-145.

\bibitem{SA1}
{\sc A. Stancu,} {\em The Discrete Planar $L_0$-Minkowski
Problem.} Adv. Math. {\bf 167} (2002),  160-174.

\bibitem{SA2}
{\sc A. Stancu}, {\em On the number of solutions to the
discrete two-dimensional $L_0$-Minkowski problem.} Adv. Math. {\bf
180} (2003), 290-323.

\bibitem{W1}
{E. Werner}, {\em Illumination bodies and affine surface area,}
Studia Math. {\bf 110}  (1994), 257-269.

\bibitem{W2}
{E. Werner}, {\em
On
$L_p$-affine surface areas}, Indiana Univ. Math. J.  {\bf 56} no. 5 (2007),  2305-2324.

\bibitem{WY} {\sc E. Werner and D. Ye}, {\em New $L_{p}$ affine isoperimetric inequalities},
Adv. Math. {\bf 218} no.3 (2008), 762-780.


\bibitem{WY1} {\sc E. Werner and  D. Ye}, {\em Inequalities for mixed $p$-affine surface area},
{Math. Ann.} {\bf  347}  (2010), 703-737



\bibitem{Z1}{\sc G. Zhang}, {\em Intersection bodies and Busemann-Petty inequalities in $\mathbb{R}^4$}, Ann.  of Math. {\bf 140} (1994), 331-346. 


\bibitem{Z2}{\sc G. Zhang}, {\em A positive answer to the Busemann-Petty problem in four dimensions}, Ann. of Math. {\bf 149} (1999), 535-543.

\bibitem{Z3}{\sc G. Zhang}, {\em New Affine  Isoperimetric Inequalities}, ICCM 2007, Vol. II, 239-267.

\vskip 1cm


\end{thebibliography}
\end{document}